\documentclass
[
11pt,
]
{amsart}
\usepackage{
amssymb
}
\newcommand{\no}[1]{#1}

\renewcommand{\no}[1]{}\newcommand{\wbar}{\bar}  \newcommand{\upDelta}{\Delta} 

\no{\usepackage{times}\usepackage[
subscriptcorrection, slantedGreek, 
nofontinfo]{mtpro}
\renewcommand{\Delta}{\upDelta}
}

\newtheorem{theorem}{Theorem}
\newtheorem{proposition}{Proposition}
\newtheorem{lemma}{Lemma}
\newtheorem{definition}{Definition}

\DeclareMathOperator{\supp}{supp}

\DeclareMathOperator{\WF}{WF_A}

\newcommand{\eps}{\varepsilon}
\newcommand{\R}{{\bf R}}
\newcommand{\Id}{\mbox{Id}}
\renewcommand{\r}[1]{(\ref{#1})}
\newcommand{\PDO}{$\Psi$DO}
\newcommand{\be}[1]{\begin{equation}\label{#1}}
\newcommand{\ee}{\end{equation}}

\DeclareMathOperator{\n}{neigh}
\renewcommand{\d}{\mathrm{d}}
\renewcommand{\i}{\mathrm{i}}

\newcommand{\bo}{\partial M}

\newcommand{\Mint}{M^\text{\rm int}}

\title[The X-ray transform for a generic family of curves]{The X-Ray transform for a generic family of curves\\ and weights}
\author[B. Frigyik]{Bela Frigyik}
\address{Department of Mathematics, Purdue University, West Lafayette, IN 47907}

\author[P. Stefanov]{Plamen Stefanov}
\address{Department of Mathematics, Purdue University, West Lafayette, IN 47907}
\thanks{Second author partly supported by NSF Grant DMS-0400869}

\author[G. Uhlmann]{Gunther Uhlmann}
\address{Department of Mathematics, University of Washington, Seattle, WA 98195}
\thanks{Third author partly supported by NSF and a Walker Family Endowed Professorship}

\begin{document}

\begin{abstract} 
We study the weighted integral transform on a compact manifold with boundary over a smooth family of curves $\Gamma$.  
We prove generic injectivity and a stability estimate under the condition that the conormal bundle  of $\Gamma$ covers $T^*M$. 
\end{abstract}

\maketitle

\section{Introduction} 
Let $M$ be a compact manifold with boundary. Let $\Gamma$ be an open family of smooth (oriented) curves on $M$, with a fixed parametrization on each one of them, with endpoints on $\bo$, such that for each $(x,\xi)\in TM\setminus 0$, there is at most one curve $\gamma_{x,\xi}\in\Gamma$ through $x$ in the direction of $\xi$, and the dependence on $(x,\xi)$ is smooth, see next section. Define the weighted ray transform
\be{01}
I_{\Gamma,w} f(\gamma) = \int w\left(\gamma(t),\dot\gamma(t)\right) f(\gamma(t))\,\d t, \quad \gamma\in\Gamma,
\ee
where $w(x,\xi)\not=0$ is a smooth non-vanishing complex valued function on $TM\setminus 0$. We study the problem of the injectivity of $I_{\Gamma,w}$ on functions on $M$. We impose no-conjugacy  conditions on $\Gamma$ that would guarantee that $I_{\Gamma,w}$ recovers singularities. Under that condition, we prove that $I_{\Gamma,w}$ is injective for generic $\Gamma$, $w$, including analytic ones, and that there is a stability estimate. 
This is a generalization of the X-ray transform arising in Computed
Tomography which consists in integrating functions over
lines provided with the standard Lebesgue measure.

In \cite{Mu}, Mukhometov showed that in a compact domain $\Omega$ in $\R^2$, $I_{\Gamma,w}$, $w=1$, is injective, for any set $\Gamma$, provided that the curves $\gamma$ have unit speed, and $\Omega$ is simple w.r.t.\ those curves. The latter means that for any two points $x$, $y$ in $\bar\Omega$, there is unique curve in $\gamma$ connecting them that depends smoothly on its endpoints. He later showed that this remains true if $w$ is close enough to a constant in an explicit way. Stability estimates were also given. In dimension $n\ge3$ there is no such known result for an arbitrary simple family of curves. On the other hand, if $\Gamma$ is the family of the geodesics of a given (simple) Riemannian or Finsler metric, and $w$ is close enough to a constant, injectivity and stability of $I_{\Gamma,w}$ was established in \cite{ Mu2, Mu3,AR, BG,R}.

The transform $I_{\Gamma,w}$ is not always injective, even for simple $\Gamma$. 
An example by Boman \cite{B} provides a smooth positive weight function $w$ so that  $I_{\Gamma,w}$ fails to be injective in a ball in $\R^2$, where $\Gamma$ consists of all straight lines.

In the present work, we have incomplete data, i.e., we do not assume that we have a curve in $\Gamma$ through any point in $M$ in the direction of any vector (unless $n=2$). On the other hand, we want $\{N^*\gamma, \; \gamma\in \Gamma\}$ to cover $T^*M$, the latter considered as a conic set. 
We do not assume convexity of the boundary w.r.t.\ $\Gamma$. If $\Gamma$ is a subset of geodesics of a certain metric, then some geodesics (not in $\Gamma$) are allowed to have conjugate points, or to be trapped, but we exclude them from $\Gamma$. On the other hand, the result is generic uniqueness and stability, and Boman's result shows that this is the optimal one in this setting.

Our approach differs from the works cited above and uses microlocal and analytic microlocal methods. Such methods are not new in integral geometry, see, e.g., \cite{Gu, GuS1, GuS2, GrU, B, BQ, Q}, but we use some recent ideas that led to new results in tensor tomography and boundary rigidity of compact Riemannian manifolds with boundary, see \cite{SU-Duke, SU-rig, SU-AJM}.

\section{Statement of the main results}
Fix a compact manifold with boundary $M_1$ such that $\Mint_1\supset M$, where $\Mint_1$ stands for the interior of $M_1$.  
We equip $M_1$ with a real analytic atlas, where $\bo$ is smooth but not necessarily analytic.  We will think of the curves $\gamma$ as extended outside $M$ to $\Mint_1$ so that their endpoints are in $\Mint_1$, and $\gamma\cap M$ remains unchanged. Different extensions will not change $I_{\Gamma,w}f$ as long as $\gamma\cap M$ is the same. 
 By $\gamma_{x,\xi}$, we will frequently denote the curve in $\Gamma$, if exists, so that $x\in\gamma_{x,\xi}$, and $\dot\gamma_{x,\xi}=\mu\xi$ at the point $x$ with some $\mu>0$. We will freely shift the parameter on $\gamma_{x,\xi}$ but not rescale it, so we may assume that $x=\gamma_{x,\xi}(0)$, then  $\dot\gamma_{x,\xi}(0)=\mu\xi$.

We want $\gamma_{x,\xi}$, for $(x,\xi)\in TM$, to depend smoothly on $(x,\xi)$, therefore in any coordinate chart, $\gamma = \gamma_{x,\xi}$ solves 
\be{11}
\ddot\gamma = G(\gamma,\dot\gamma),
\ee
where $G(x,\xi)= \ddot\gamma_{x,\xi}(0)$ is smooth. The generator $G(x,\xi)$ is only defined for $|\xi|=|\dot\gamma_{x,\xi}(0)|$ (in any fixed coordinates) but we can extend it for all $\xi$. In case of a Riemannian metric, for example, $G^i(x,\xi) = -\Gamma^i_{kl}(x)\xi^k\xi^l$, for $|\xi|_g=1$, and extended for all $\xi$. The generator $G$ determines a vector field $\mathbf{G}$ on $TM$ that in local coordinates is given by 
\be{12}
\mathbf{G} = \xi^i\frac{\partial}{\partial x^i} + G^i(x,\xi) \frac{\partial}{\partial \xi^i},
\ee
see also \r{41}, \r{42}. The curves $\gamma\in\Gamma$ are the projections of integral curves of $\mathbf{G}$ to the base, with appropriate initial conditions that reflect the choice of the parametrization.

We assume that $\Gamma$ is open with a natural smooth structure as follows. Fix any $\{\gamma(t);\; l^-\le t\le l^+\}\in\Gamma$, $\pm l_\pm>0$, 
$\gamma(l_\pm)\in\Mint_1\setminus M$, where we shifted the parameter $t$ arbitrarily, and set $x_0:=\gamma(0)$. Let  $H$ be a hypersurface in $\Mint_1\setminus M$ intersecting $\gamma$ transversally at $x_0$, and let $\xi_0=\dot\gamma(0)$. We assume that  there exists a neighborhood $U$ of $(x_0,\xi_0)$ and a smooth positive function  $\mu(x,\xi)$, 
$(x,\xi)\in U\cap H$, with $\mu(x_0,\xi_0)=1$, 
so that the integral curves of $\mathbf{G}$ with initial conditions $(\gamma(0), \dot \gamma(0) ) = (x,\mu(x,\xi)\xi)$, $(x,\xi)\in U\cap H$, and interval of definition $l^-\le t\le l^+$ belong to $\Gamma$ (and in particular, the endpoints are in $\Mint_1\setminus M$). 
This makes $\Gamma$ a smooth manifold; if $H$ is given locally by $x^n=0$, then 
$\Gamma$ is locally parametrized by $(x',\theta)\in \R^{n-1}\times S^{n-1}$. We say that $\Gamma$ is $C^k$, respectively \textit{analytic}, if $G$ is $C^k$,
respectively analytic, on $TM_1$, and for any such choice of $C^k$, respectively analytic $H$, the functions $\mu$ are $C^k$, respectively analytic, too. 

It is not hard to see that by duality, one can define $I_{\Gamma,w}f$ for any distribution $f\in \mathcal{D}'(\Mint_1)$ supported in $M$. 

Given $x\in M$, we define the exponential map $\exp_x(t,\xi)$, $t\in \R$, $\xi\in TM\setminus 0$, as $\exp_x(t,\xi)=\gamma_{x,\xi}(t)$. Note that $\exp_x(t,\xi)$ is a positively homogeneous function of order 0 in the $\xi$ variable, and in local coordinates, we can think that $\xi\in S^{n-1}$. Then $x=\gamma(0)$ and $y=\gamma(t_0)$ will be called \textit{conjugate} along $\gamma$, if $D_{t,\xi}\exp_x(t,\xi)$ has rank less than $n$ at $(t_0,\xi_0)$, where $\xi_0=\dot\gamma(0)$. It is easy to see that this definition is independent of a change of the parametrization along the curves in $\Gamma$ (that we keep fixed). We would like to note here that (in a fixed coordinate system), the map $v=t\xi\mapsto \exp_x(t,\xi)$, where $|\xi|=1$, $t\in \R$, may not be $C^\infty$. In case of magnetic systems, for example, it is only $C^1$ while $\exp_x(t,\xi)$ is a smooth function of all variables, see \cite{DPSU}. This requires some modifications in the analysis of the normal operator \r{Na}, see section~\ref{singular}.

It is clear that one cannot hope to recover any $f$ from $I_{\Gamma,w}f$, if there is a point in $M$ so that no $\gamma\in\Gamma$ goes through it. We impose a microlocal condition that  requires something more than that, we want any $(x,\zeta)\in T^*M\setminus 0$ to be ``seen'' by some simple $\gamma\in\Gamma$. 

\begin{definition}   \label{def_1}
We say that $\Gamma$ satisfying the assumptions above is a \textbf{regular} family of curves, if for any $(x,\zeta)\in T^*M$, there exists $\gamma\in\Gamma$ through $x$ normal to $\zeta$ without conjugate points.

We call any $\gamma$ as above a \textbf{simple} curve. 
\end{definition} 

If $\Gamma$ is not regular, one can give the following example of a non-injective $I_{\Gamma,w}$. Let $M$ be a subdomain with boundary of the sphere $S^{n-1}$ with its natural metric. Assume that $\Mint$ contains a pair of antipodal points $a$ and $b$. Then any function that is supported in two symmetric to each other small enough neighborhoods $A\ni a$, $B\ni b$, and odd with respect to the antipodal map, integrates to 0 over any geodesic in $M$. Not only $I_{\Gamma,w}f$ with $w=1$ does not determine $f$, it does not determine the singularities, either. For example, if $f=\delta_a-\delta_b$, where $\delta_{a,b}$ are delta distributions centered at $a$ and $b$, respectively; then $I_{\Gamma,1}f=0$. 

On the other hand, one can see that $I_{\Gamma,w}f$, known for a regular family of curves resolves the singularities of $f$. Using analytic microlocal arguments, we also show that one can recover the analytic singularities, as well, if $\Gamma$ is analytic. This allows us to prove the following.

\begin{theorem} \label{thm_an}
Let $\Gamma$ be an analytic regular family of curves in $M_1$ and let $w$ be analytic and non-vanishing in $M$. Then $I_{\Gamma,w}f=0$ for   $f\in \mathcal{D}'(\Mint_1)$ supported in $M$ implies $f=0$. In particular, $I_{\Gamma,w}$ is injective on $L^1(M)$.
\end{theorem}

To formulate a stability result, we will fix a parametrization of $\Gamma$. Let $H$ be a finite collection of hypersurfaces $\{H_m\}$ in $\Mint_1$ that are allowed to intersect each other. Then $H$ may not be a hypersurface but is still a manifold if we think of each $H_m$ as belonging to a different copy of $M$. 
Let $\mathcal{H}$ be an open conic subset of $\{(z,\theta)\in TM_1; \;  z\in H,  \;  \theta\not\in T_zH\}$, and let $\pm l^\pm(z,\theta)\ge0$ be two continuous functions. Let $\Gamma(\mathcal{H})$ be the subset of curves  of $\Gamma$  originating from $\mathcal{H}$, i.e., 
\be{5}
\Gamma(\mathcal{H}) = \left\{\gamma_{z,\theta}(t); \; l^-(z,\theta)\le t\le l^+(z,\theta), \; (z,\theta)\in \mathcal{H}  \right\}.
\ee
We also assume that each $\gamma\in \Gamma(\mathcal{H})$ is a simple curve. 

We will fix a parametrization of a subset of $\Gamma$ that is still regular. 

Given $\mathcal{H}$ as above, we consider an open set $\mathcal{H}' \Subset \mathcal{H}$, and let $\Gamma(\mathcal{H}')\Subset \Gamma(\mathcal{H})$ be the associated set of  curves defined as in \r{5}, with the same $l^\pm$.  The restriction $\gamma\in \Gamma(\mathcal{H}')\subset \Gamma(\mathcal{H})$ can be modeled by introducing a weight function $\alpha$ in $\mathcal{H}$, such that $\alpha=1$ on $\mathcal{H}'$, and $\alpha=0$ otherwise. It is more convenient to allow $\alpha$ to be smooth but still supported in $\mathcal{H}$. 

We consider $I_{\Gamma,w,\alpha}=\alpha I_{\Gamma,w}$, or more precisely, 
\be{I_a0}
I_{\Gamma,w,\alpha}f = \alpha(z,\theta) \int_{l^-(z,\theta))}^{l^+(z,\theta)} w\big(\gamma_{z,\theta}, \dot \gamma_{z,\theta}\big) f(\gamma_{z,\theta}) \,\d t, \quad (z,\theta)\in  \mathcal{H}.
\ee
Next, we set
\be{Na}
N_{\Gamma,w,\alpha} = I_{\Gamma,w,\alpha}^*I_{\Gamma,w,\alpha} 
= I_{\Gamma,w}^*|\alpha|^2I_{\Gamma,w}.
\ee
Here the adjoint is taken w.r.t.\ a fixed positive smooth measure $\d\Sigma$ on $\mathcal{H}$; more precisely, we assume that in any local coordinate chart, $\d\Sigma := \sigma(z,\theta) \,\d S_z\,\d \theta$ on $\mathcal{H}$,  where $\d S_z$ is the surface measure on $H$ in the so fixed coordinate system,  $\d \theta$ is the surface measure on $S^{n-1}$, and $C^\infty\ni\sigma>0$. Notice that $\d\Sigma$ is not invariant under a different choice of $\mathcal{H}$ and a coordinate system on it. On the other hand, injectivity of $N_{\Gamma,w,\alpha}$ is equivalent to injectivity of $I_{\Gamma,w,\alpha}$, and the latter is equivalent to injectivity of $I_{\Gamma,w}$ restricted to $\supp\alpha$, see \cite{SU-Duke}, and this property is independent of the choice of $\mathcal{H}$ and the coordinates on it as long as they parametrize the same set of curves.

\begin{theorem}   \label{thm_2} \ 

(a) %
Let $\mathcal{H}'\Subset\mathcal{H}$ be as above with $\Gamma(\mathcal{H}')\subset \Gamma(\mathcal{H})$ regular, and $(G,\mu,\sigma,w)$ fixed.   Fix $\alpha\in C^\infty$ with  $\mathcal{H}'  \subset\supp\alpha\subset \mathcal{H}$. 
If $I_{\Gamma,w,\alpha}$ is injective, where $\Gamma =\Gamma(\mathcal{H})$, then we have 
\be{est}
\|f\|_{L^2(M)} /C   \le 
\|N_{\Gamma,w,\alpha} f\|_{ H^1(M_1)} \le C\|f\|_{L^2(M)}.
\ee

(b) Let $\mathcal{H}'\Subset\mathcal{H}$, $\alpha= \alpha^0$ be as above related to some fixed $(G_0,\mu_0,\sigma_0,w_0)$. Assume that $I_{\Gamma_{0},w_0,\alpha^0}$ is injective, where $\Gamma_0 =\Gamma_0(\mathcal{H})$. Then 
 estimate \r{est} remains true for  $(G,\mu,\sigma, w,\alpha)$  belonging to a small $C^2$ neighborhood of $(G_0,\mu_0,\sigma_0,w_0,\alpha^0)$,  with a uniform constant $C>0$. 
\end{theorem}

\noindent{\bf Remark}
In fact we need only $C^1$ regularity for $w$, $\alpha$.

We notice that $C^2$ above refers to different spaces. More precisely, $\mu$, $\alpha^0$ are considered in $C^2(\mathcal{H})$, while $G$, $w$ are considered in $C^2(TM)$. To define correctly $C^2(TM)$, we fix any finite atlas on $M$, see also the remark in section~\ref{sec_4}. 

\medskip
\textbf{Example (simple systems).} Let $M\subset \R^n$ be diffeomorphic to a ball, and let $G(x,\xi)$ be a smooth generator on $TM\setminus 0  \cong M\times \R^n\setminus 0$. 
Fix a coordinate system on $M$. We can assume that $G$ is defined on $SM \cong M\times S^{n-1}$ and extend as a homogeneous of order 0 to all $\xi\not=0$. Set 
\[
\partial_-SM = \left\{ z\in \bo;\; \theta\cdot\nu<0  \right\}
\]
where $\nu(z)$ is the exterior unit normal to $\bo$.  Then we define $\Gamma$ as the set of all curves  $\gamma = \gamma_{z,\theta}$  that solve
\be{26}
\ddot\gamma = G(\gamma,\dot\gamma), \quad \gamma(0)=z, \quad \dot\gamma(0) = \lambda(z,\theta)\theta,  \quad (z,\theta)\in \partial_-SM, 
\ee
where $\lambda>0$ is a given smooth function on $M\times S^{n-1}$ with $\lambda(z,\theta) = |\dot\gamma_{z,\theta}(0)|$. Let 
$\gamma_{z,\theta}$ be the maximal curves with those initial conditions.  Assume that for any $x\in M$, the map $\exp_x : \exp_x^{-1}(M) \to M$ is a diffeomorphism depending smoothly on $x$. Note that this implies that all those curves are of finite length;  for any $x$, $y$ in $M$, there is unique $\gamma\in \Gamma$ that passes through them, smoothly depending on $x$, $y$, and the curves in $\Gamma$ have no conjugate points. As above, $\gamma$'s are allowed to be directed curves; if $x\in \Mint$, $\theta\in S^{n-1}$ then the curves $\gamma_{x,\theta}$ and $\gamma_{x,-\theta}$ are not necessarily the same. We also assume that $M_1\Supset M$ (meaning that $\Mint_1\supset \bar M=M$) is another domain diffeomorphic to a ball so that $(G,\lambda)$ extends smoothly there and satisfies the same assumptions. 

For a simple system as above, define
\be{27}
I_{G,\lambda,w}f(z,\theta) = \int  w\left(\gamma_{z,\theta}(t),\dot\gamma_{z,\theta}(t)\right) f(\gamma_{z,\theta}(t))\,\d t, \quad (z,\theta) \in  \partial_-SM_1.
\ee
One could also  study subsets of curves as above. 
Let $\sigma$ be any positive $C^1$ function on $\overline{\partial_-SM_1}$, and set $\d\Sigma = \sigma(z,\theta)|\nu(z)\cdot\theta|\,\d S_z\,\d\theta$. Then
\[
I_{G,\lambda,w} : L^2(M) \to L^2(\partial_-SM_1,\d\Sigma)
\]
is a bounded map,  and $N_{G,\lambda,w} = I_{G,\lambda,w}^*I_{G,\lambda,w}$ is a well defined operator on $L^2(M)$ that can be extended as an  operator from $L^2(M)$ to $H^1(M_1)$. Note that the factor $|\nu(z)\cdot\theta|$ in $\d \Sigma$ can be omitted since $\bo_1$ is convex and  $M$ stays at a positive distance from $\bo_1$. If $M_1=M$, and if $\bo$ is strictly convex w.r.t.\ $\Gamma$, then that factor is needed to preserve the mapping properties of $N_{G,\lambda,w}$; see \cite{SU-Duke} for the Riemannian case.


\section{Injectivity of $I_{\Gamma,w}$ for analytic systems}
In this section we prove Theorem~\ref{thm_an}. We denote by $\WF(f)$ the analytic wave front set of $f$. 

\begin{proposition}  \label{prop_wf}
Let $\gamma_0\in\Gamma$ be a simple curve. Let $I_{\Gamma,w}f(\gamma)=0$ for some $f\in\mathcal{D}'(M_1)$ with $\supp f\subset M$ and all $\gamma\in\n(\gamma_0)$. Let $\Gamma$ and $w$ be analytic near $\gamma_0$. Then
\be{31a}
N^*\gamma_0 \cap \WF(f) =\emptyset.
\ee
\end{proposition}

\begin{proof} 
We will choose first a coordinate system $(x',x^n)$ near $\gamma_0$ so that the latter is given by $x'=0$, $x^n=t$, $t\in [l^-,l^+]$ with some $\pm l^\pm\ge0$, and moreover, replacing $x'=0$ by $x'=z$, where $z$ is a constant vector with $|z|\ll1$, one still gets a curve in $\Gamma$ (parametrized by $t$ again, i.e., a unit speed line segment in the so fixed coordinate system). 

Fix a point $p_0\in\gamma_0$, and shift the parametrization of $[l^-,l^+]\ni t\mapsto \gamma_0(t)$ so that $p_0=\gamma_0(0)$. Assume that $p_0\not\in M$ and  that the part of $\gamma_0$ corresponding to $l^-\le t\le 0$ is outside $M$, too. Set $x=\exp_{p_0}(t,\theta)$, where $|\theta|=1$, $t\ge0$, where the norm is in any  fixed coordinate system near $p_0$.  Then $(t,\theta)$ are local coordinates near any point on $\gamma_0\cap M$ because the $\gamma_0$ is simple. Since $\gamma_0$ may self-intersect, they may not be global ones. On the other hand, there can be finitely many intersections only, and one can assume that each time $\gamma_0$ intersect itself, it happens on a different copy of $M\times\R$. More precisely, $(t,\theta) \mapsto (t,\exp_{p_0}(t,\theta)) \subset M\times\R$ is a codimension one submanifold of $M$ for $\theta$ close to $\theta_0=\dot\gamma_0(0)$ and  $t\in (0,l^+)$  by the simplicity assumption, and we think of any function $f:M\to \mathbf{C}$ as defined on that manifold. 
Therefore, without loss of generality, we may assume that  $\gamma_0$ does not self-intersect. 

Write $x'=\theta'$, $x^n = t$. 
Then $x$ are the coordinates we were looking for in
\[
U = \left\{ x; \;  |x'|<\eps, \, l^-<t<l^+ \right\} \subset M_1
\]
with $0<\eps\ll1$. They are analytic, since $\Gamma$ is analytic. 

Fix $x_0\in\gamma_0$, and $\xi^0\in T^*M_1$ conormal to $\gamma_0$. We need to prove that 
\be{31}
(x_0,\xi^0) \not\in \WF(f).
\ee
By shifting the $x^n$ coordinate, we can always assume that $x_0=0$. Note that $\theta_0 :=\dot\gamma_0(0)= e_n$. Here and below, $e_j$ stand for the vectors $\partial/\partial x^j$, and $e^j$ stand for the covectors $dx^j$.

Assume first that $f$ is  continuous in $M$ and vanishes outside $M$. 

The arguments that follow are close to those in \cite{SU-AJM}. 
Set first $Z = \{x^n=0;\; |x'|<7\varepsilon/8\}$, and denote the $x'$ variable on $Z$ by $z'$. We will work with the curves $t\mapsto \gamma_{(z',0),(\theta',1)}(t)$ 
defined on $l^-\le t\le l^+$, the same interval on which $\gamma_0$ is defined. 
Each such  curve is in $\Gamma$ for $|\theta'|\ll1$ because the latter is open. They all have endpoints in $\Mint_1\setminus M$, and in fact, we modified a bit the endpoints of the interval of definition to make them constant ($l^\pm$). We can do this, when $\eps\ll1$, and this does not affect integrals of $f$ over them.

Let $\chi_N(z')$, $N=1,2,\dots$,  be a sequence of smooth cut-off functions equal to $1$ for $|z'|\le 3\varepsilon/4$, supported in $Z$, and satisfying the estimates 
\be{N}
\left| \partial^\alpha \chi_N\right|\le (CN)^{|\alpha|}, \quad 
|\alpha|\le N,
\ee
see \cite[Lemma~1.1]{T}. 
Set $\theta=(\theta',1)$, $|\theta'|\ll1$,  and multiply 
$$
I_{\Gamma,w} f \left(\gamma_{(z',0),\theta}\right) =0
$$
by $\chi_N(z') e^{\mathrm{i}\lambda z'\cdot \xi'}$, where $\lambda>0$, $\xi'$ is in a complex neighborhood of $(\xi^0)'$, and integrate w.r.t.\ $z'$ to get
\be{42'}
\iint e^{\mathrm{i}\lambda   z'\cdot \xi'} \chi_N(z') w\left(\gamma_{(z',0),\theta}(t), \dot\gamma_{(z',0),\theta}(t) \right) f\left( \gamma_{(z',0),\theta} (t)\right) (t)\, \d t\, \d z'=0.
\ee

For $|\theta'|\ll1$,  $(z',t)\in Z\times (l^-,l^+)$ are local coordinates near $\gamma_0$ given by $x=\gamma_{(z',0),\theta} (t)$. Indeed, if $\theta'=0$, we have $x=(z',t)$. Therefore, for $\theta'$ fixed and small enough, $(t,z')$ are analytic local coordinates, depending analytically on $\theta'$.  In particular,  $x=(z'+t\theta',t) + O(|\theta'|)$. Performing a change of variables in \r{42'}, we get
\be{43}
\int e^{\mathrm{i}\lambda z'(x,\theta')\cdot \xi'}  a_N(x,\theta') f( x) \, \d x=0
\ee
for $|\theta'|\ll1$, $\forall\lambda$, $\forall\xi'$, where, for $|\theta'|\ll1$,  the function $(x,\theta') \mapsto a_N$ is analytic, independent of $N$, and non-zero for $x$ in a neighborhood of $\gamma_0$, satisfies \r{N} everywhere, vanishes for  $x\not\in U$; 
and $a_N(0,\theta')=w(0,\theta)$. 

Without loss of generality we can assume that 
\[
\xi^0 =e^{n-1}.
\]
Here and below, $e_j$ stand for the vectors $\partial/\partial x^j$, and $e^j$ stand for the covectors $dx^j$. 

We choose the following vector $\theta(\xi)$  analytically depending on $\xi$ near $\xi=\xi^0$:
\be{45_1}
\theta(\xi) = \bigg( \xi_1,\dots,\xi_{n-2}, -\frac{\xi_1^2+\dots+\xi_{n-2}^2 +\xi_n}{\xi_{n-1}} ,1    \bigg).
\ee
If $n=2$, this reduces to $\theta(\xi) = (-\xi_2/\xi_1,1)$.  
Clearly, 
\be{th}
\theta(\xi)\cdot\xi=0, \quad \theta^n(\xi)=1, \quad \theta(\xi^0) = e_n.
\ee

Differentiate \r{45_1} to get
\be{45_2}
\frac{\partial\theta}{\partial \xi_\nu}(\xi^0) = e_\nu, \quad \nu=1,\dots,n-2, \quad \frac{\partial\theta}{\partial \xi_{n-1}}(\xi^0)=0,  \quad   \frac{\partial\theta}{\partial \xi_{n}}(\xi^0) = -e_{n-1}.
\ee
In particular, the differential of the map $S^{n-1}\ni \xi \mapsto \theta'(\xi)$ is invertible at $\xi=\xi^0=e^{n-1}$.

Replace $\theta=(\theta',1)$ in \r{43} by $\theta(\xi)$  (the requirement $|\theta'| \ll1$ is fulfilled for $\xi$ close enough to $\xi^0$), to get
\be{45}
\int e^{\mathrm{i}\lambda \varphi(x,\xi)} \tilde a_N(x,\xi)  f(x) \, \d x=0,  \
\ee
where $\varphi$ is analytic in $U$, and 
$\tilde a_N$ has the properties of $a_N$ above for $\xi$ close enough to $\xi^0$.  
In particular,
$$
\tilde a_N(0,\xi)=w(0,\theta(\xi)). 
$$
The phase function is given by 
\be{45a}
\varphi(x,\xi) = z'(x,\theta'(\xi))\cdot \xi'.
\ee
To verify that $\varphi$ is a non-degenerate phase in  $\n(0,\xi^0)$, 
i.e., that $\det \varphi_{x\xi}(0,\xi^0)\not =0$, note first that $z'=x'$ when $x^n=0$, therefore, $(\partial z'/\partial x')(0,\theta(\xi))=\Id$. On the other hand, linearizing near $x^n=0$, we easily get $(\partial z'/\partial x^n)(0,\theta(\xi))=-\theta'(\xi)$. Therefore,
\be{45b}
\varphi_x(0,\xi) = (\xi', -\theta'(\xi)\cdot \xi') = \xi
\ee
by \r{th}. So we get $\varphi_{x\xi}(0,\xi) = \Id$, which proves the non-degeneracy claim above. In particular,  $x\mapsto \varphi_\xi(x,\xi)$ is a local diffeomorphism in $\n(0)$ for $\xi\in\n(\xi^0)$, and therefore injective. We need however a semiglobal version of this along $\gamma_0$ as in the lemma below.

\begin{lemma}  \label{lemma_phase} 
There exists $\delta>0$ such that 
\[
\varphi_\xi(x,\xi) \not= \varphi_\xi(y,\xi) \quad \text{for $x\not=y$},
\]
for   $x\in U$, $|y|<\delta$, $|\xi-\xi^0|<\delta$, $\xi$ complex.
\end{lemma}

\begin{proof}
We will prove the lemma first for $y=0$, $\xi=\xi^0$, $x'=0$. Since $\varphi_\xi(0,\xi)=0$, we need to prove that the only solution to $\varphi_\xi((0,x^n),\xi^0)=0$ in the interval $l^- \le x^n \le l^+$ is $x^n=0$.

We start with the observation that $\varphi(\gamma_{0,(\theta'(\xi),1)}(t),\xi)=0$. Differentiate the latter w.r.t.\ $\xi$ at $\xi=\xi^0$, $t=x^n$, to get
\[
\frac{\partial\varphi}{\partial\xi_i}((0,x^n),\xi^0) = -\frac{\partial}{\partial \xi_i}\bigg|_{\xi=\xi^0}   \varphi\big(\gamma_{0,(\theta'(\xi),1)}(x^n),\xi^0\big)
= -\frac{\partial\varphi}{\partial x^j}((0,x^n),\xi^0 )J^j_\nu(0,x^n) \frac{\partial\theta^\nu}{\partial\xi_i}(\xi^0),
\]
where $J_\nu(t) = \partial \gamma_{0,\theta}(t)/\partial\theta_\nu$ at $\theta=e_n$, $\nu=1,\dots,n-1$, are ``Jacobi'' vector fields. Since $\varphi(x,\xi^0) = x'\cdot(\xi^0)'=x^{n-1}$, we get by \r{45_2}, (recall that $\xi^0=e^{n-1}$),
\be{45_3}
\frac{\partial\varphi}{\partial\xi_j} ((0,x^n),\xi^0) =
\begin{cases}
-J_j^{n-1}(x^n), & j =1,\dots,n-2,\\
                        0,  & j=n-1,\\
J_{n-1}^{n-1}(x^n),  & j =n,
\end{cases}
\ee
where $J^{n-1}_\nu$ is the $(n-1)$-th component of $J_\nu$. 
Now, assuming that the l.h.s.\ of \r{45_3} vanishes for some fixed $x^n=t_0$, we get that $J_\nu^{n-1}(t_0)=0$, $\nu=1,\dots,n-1$. On the other hand, $\Sigma := \text{span}(J_1(t_0), \dots, J_{n-1}(t_0))$ is a hyperplane transversal to $e_n$ by the simplicity assumption. Therefore, for the unit normal $\nu$ to $\Sigma$, we have $\nu_n\not=0$. Hence, $\nu$ and $e_{n-1}$ are linearly independent, and the intersection of $\Sigma$ and $e_{n-1}^\perp$ is of codimension 2, and $J_1(t_0), \dots, J_{n-1}(t_0)$ all belong there.  Therefore, $J_\nu(t_0)$, $\nu=1,\dots,n-1$, form a linearly dependent system of vectors. 
The latter contradicts the simplicity assumption.

The same proof applies if $x'\not=0$ by shifting the $x'$ coordinates. 

Let now $y$, $\xi$ and $x$ be as in the Lemma. The lemma is clearly true for $x$  in the ball $B(0,\eps_1) = \{|x|<\eps_1\}$, where $\eps_1\ll1$, because $\varphi(0,\xi^0)$ is non-degenerate. On the other hand, $\varphi_\xi(x,\xi)\not= \varphi_\xi(y,\xi)$ for $x\in \bar U\setminus B(0,\eps_1)$,  $y=0$, $\xi=\xi^0$. Hence, we still have $\varphi_\xi(x,\xi)\not= \varphi_\xi(y,\xi)$ for a small perturbation of  $y$ and $\xi$.
\end{proof}

We will apply the complex stationary phase method \cite{Sj-Ast}, see also \cite[Section~6]{KSU}. For $x$, $y$ as in Lemma~\ref{lemma_phase}, and   $|\eta-\xi^0|\le \delta/\tilde C$, $\tilde C\gg1$, $\delta\ll1$, multiply \r{45} by 
$$
\tilde \chi(\xi-\eta)e^{\i \lambda( \i (\xi-\eta)^2/2  -\varphi(y,\xi) )},
$$
where $\tilde \chi$ is the characteristic function the complex ball $B(0,\delta)$, and integrate w.r.t.\ $\xi$ to get
\be{45_4}
\iint e^{\i\lambda \Phi(y,x,\eta,\xi)} b_N(x,\xi,\eta)  f( x) \, \d x\, \d \xi=0, 
\ee
where $b_N$ is another amplitude, analytic, independent of $N$,   and elliptic near $\gamma_0\times \{\xi^0\}$, satisfying \r{N}, and
\[
\Phi = -\varphi(y,\xi)+\varphi(x,\xi) +\frac{\i}2 (\xi-\eta)^2.
\]
We study the critical points of $\xi\mapsto \Phi$.  If $y=x$, there is a unique (real) critical point $\xi_{\rm c}=\eta$, and it satisfies  $\Im\Phi_{\xi\xi} >0$ at $\xi= \xi_{\rm c}$. For $y\not=x$, there is no real critical point by Lemma~\ref{lemma_phase}. On the other hand, again by Lemma~\ref{lemma_phase}, there is no (complex) critical point if $|x-y|>\delta/C_1$ with some $C_1>0$, and there is 
a unique complex critical point $\xi_{\rm c}$  if $|x-y|<\delta/C_2$, with some $C_2>C_1$, still non-degenerate if $\delta\ll1$. For any $C_0>0$, if we integrate in \r{45_4} for $|x-y|>\delta/C_0$, and use the fact that $|\Phi_\xi|$ has a positive lower bound (for $\xi$ real), we get
\be{45_5}
\bigg| \iint_{|x-y|>\delta/C_0} e^{\i\lambda \Phi(y,x,\eta,\xi)} b_N(x,\xi,\eta)  f( x) \, \d x\, \d \xi \bigg| \le  C_3(C_3N/\lambda)^N  +CNe^{-\lambda/C}.
\ee
Estimate \r{45_5} is obtained by integrating $N$ times by parts, using the identity
\[
Le^{\i\lambda \Phi} = e^{\i\lambda \Phi}, \quad L := \frac{\wbar\Phi_\xi\cdot \partial_\xi}{\i\lambda|\Phi_\xi|^2}
\]
as well as using the estimate \r{N}, and the fact that on the boundary of integration in $\xi$, the $e^{\i\lambda\Phi}$ is exponentially small. 
Choose $C_0\gg C_2$. Note that  $\Im \Phi>0$ for  $\xi\in\partial (\supp \tilde\chi(\cdot-\eta))$, and $\eta$ as above, as long as $\tilde C\gg1$, and by choosing $C_0\gg1$, we can make sure that $\xi_{\rm c} $ is as close to $\eta$, as we want. 

To estimate \r{45_4} for $|x-y|<\delta/C_0$, set 
$$
\psi(x,y,\eta) := \Phi\big|_{\xi=\xi_{\rm c}}.
$$
Note that $\xi_{\text{c}} =-\i(y-x)+\eta+O(\delta)$, and $\psi(x,y,\eta) = \eta\cdot (x-y) +\frac{\i}2 |x-y|^2+O(\delta)$. The stationary complex phase method \cite{Sj-Ast}, see Theorem~2.8 there and the remark after it,  together with \r{45_5}, gives
\be{47}
\int_{|x-y|\le \delta/C_0}  e^{\i \lambda \psi(x,\alpha)}  
  f( x) B(x,\alpha; \lambda)  \, \d x =  O\big( Ne^{-\lambda/C} + \lambda^{n/2} (C_3N/\lambda)^N \big),  \quad \forall N,
\ee 
where $\alpha = (y,\eta)$, and $B$ is a classical elliptic analytic symbol \cite{Sj-Ast}, independent of $N$. Moreover, the principal symbol  $\sigma_p(B)(0,0,\eta)$ equals $w(0,\theta(\eta))$ times an elliptic factor, and is therefore elliptic itself.   Recall that $w(0,\theta(\xi^0)) =w(0, e_n)\not=0$. Take $N$ so that $N \le\lambda/(C_3e)\le N+1$ to conclude that the r.h.s.\ of \r{47} is $O(e^{-\lambda/C})$. 

At $y=x$ we have
\be{46a}
\psi_y(x,x,\eta) = -\varphi_x(x,\eta), \quad
\psi_x(x,x,\eta) =  \varphi_x(x,\eta), \quad
\psi(x,x,\eta)=0.
\ee
We also get that 
\be{46b}
\Im \psi(y,x,\eta) \ge |x-y|^2/C,
\ee
that can be obtained by writing $y=x+h$, and expanding $\psi$ in terms of powers of $h$ up to $O(h^3)$. 

Define the transform
$$
\alpha \longmapsto \beta = \left(\alpha_x, \nabla_{\alpha_x}\varphi(\alpha)\right),
$$
where, following \cite{Sj-Ast}, $\alpha=(\alpha_x,\alpha_\xi)$. 
This is equivalent to setting $\alpha=(y,\eta)$, $\beta = (y,\zeta)$, where  $\zeta = \varphi_y(y,\eta)$. Note that $\zeta =\eta+O(\delta)$, and at $y=0$, we have $\zeta=\eta$, by \r{45b}. 
It is a diffeomorphism from  a neighborhood of $(0,\xi^0)$    to its image, leaving $(0,\xi^0)$ fixed. Denote the inverse map by $\alpha(\beta)$. Note that this map and its inverse  preserve the first (n-dimensional) component and change only the second one.  
Plug $\alpha=\alpha(\beta)$ in \r{47} to get 
\be{48}
\int_{|x-\alpha_x|\le \delta/C_0}  e^{\i \lambda \psi(x,\beta)} 
B(x,\beta; \lambda) f( x) \, \d x =  O\big( e^{-\lambda/C} \big),  
\ee
for $\beta\in \n(0,\xi^0)$, where $\psi$,  $B$ are (different) functions having the same properties as above, except that now $\psi$ satisfies 
\be{49}
\psi_y(x,x,\zeta) = -\zeta, \quad
\psi_x(x,x,\zeta) =  \zeta, \quad
\psi(x,x,\zeta)=0.
\ee
By \cite[Definition~6.1]{Sj-Ast}, \r{46b}, \r{48}, \r{49}, together with the ellipticity of $B$ imply that 
\[
(0,\xi^0)\not\in \WF(f).
\]
Note that in \cite{Sj-Ast}, it is  required that $f$ must be replaced by $\bar f$ in \r{48}. If $f$ is complex-valued, we could use the fact that $I(\Re f)(\gamma)=0$, and $I(\Im f)(\gamma)=0$ for $\gamma$ near $\gamma_0$ and then work with real-valued $f$'s only. 

If $f$ is a distribution, then one can see that \r{43} still remains true with the integral in the $x$ variable understood in distribution sense. The rest of the proof remains the same, except that the cutoffs w.r.t.\ $x$ have to be replaced by smooth ones. The characterization of $\WF(f)$ in \cite[Definition~6.1]{Sj-Ast}  is formulated for distributions, too.

This concludes the proof of Proposition~\ref{prop_wf}. 
\end{proof}

\begin{proof}[Proof of Theorem~\ref{thm_an}]  
The proof of Theorem~\ref{thm_an} now follows immediately. By Proposition~\ref{prop_wf}, $f$ is analytic in $M_1$ and has compact support there. Therefore, $f=0$. 
\end{proof}

\section{The smooth parametrix}  \label{sec_4}

Under coordinate changes $x\mapsto x'$, $\mathbf{G}$ preserves its form, i.e., 
\be{41}
\mathbf{G} = \xi'^i \frac{\partial}{\partial x'^i} + G'^i(x',\xi')\frac{\partial}{\partial\xi'^i} \,,
\ee
and the transformation law is
\be{42}
G'^i(x',\xi') = G^k\Big(x',\frac{\partial x}{\partial x'^j}\xi'^j\Big) \frac{\partial x'^i}{\partial x^k} 
+ \frac{\partial^2 x'^j}{\partial x^i\partial x^k} \frac{\partial x^i}{\partial x^s} \frac{\partial x^k}{\partial x^t} \xi'^s\xi'^t.
\ee
This shows that the assumption $G\in C^k$ is independent of the choice of the coordinate chart, and choosing a different finite atlas will preserve inequalities of the kind $\|G-\tilde G\|_{C^2(TM)} \le C\eps$ by changing $C$ only. 
\medskip

We construct below a parametrix for $N_{\Gamma,w,\alpha}$ assuming that $G, \lambda, w,\alpha$ are smooth. 

\begin{proposition}  \label{pr_2}
$N_{\Gamma,w,\alpha}$ is an elliptic  classical \PDO\ of order $-1$ in $\Mint$. As a consequence, there exists a classical pseudodifferential operator $Q$ in $\Mint_1$ of order $1$ so that
\[
QN_{\Gamma,w,\alpha} f = f+Kf
\]
for any $f\in \mathcal{D}'(\Mint_1)$ with $\supp f\subset M$, and an operator $K$ with a $C_0^\infty(\Mint_1\times\Mint_1)$ Schwartz kernel. 
\end{proposition}

As a first step towards the proof of Proposition~\ref{pr_2},  we derive a formula for $I^*_{\Gamma,w,\alpha}$. 
Notice that the map $\mathcal{H}\times (l^-,l^+)\ni(z,\theta,t)\mapsto (x,v)\in \R^n\times S^{n-1}$ given by $x= \exp_z(t,\theta)$, $v = \partial_t   \exp_z(t,\theta)/|\partial_t \exp_z(t,\theta) |$ is a local diffeomorphism. Indeed, fix $(z_0, \theta_0, t_0)$, and let $(x_0,v_0)$ be the corresponding $(x,v)$. To find the inverse of that map, we need to solve 
\[
\exp_x(-t,v) = z, \quad -\partial_t \exp_x(-t,v) = \theta, \quad z\in H,
\]
for $(z,\theta,t)$ near $(z_0, \theta_0, t_0)$, so that $(z,\theta,t)= (z_0, \theta_0, t_0)$ for $(x,v)=(x_0,v_0)$. This can be done, since $H$ is not tangent to any $\theta$  such that $(z,\theta)\in \mathcal{H}$. Let  $J^\flat(x,v) = \d(z,\theta,t)/ \d (x,v)$ be the corresponding Jacobian (depending on the choice of the local chart near $x$).

Let $\phi\in C_0^\infty(\mathcal{H})$,  $f\in C^\infty(M)$, and let $w_1\in C^\infty_0(T\Mint_1)$ have small enough support that can fit in a coordinate chart that we fix. Then
\[
\begin{split}
\int (I_{\Gamma,w_1,\alpha} f) \bar\phi\, \d\Sigma &= \iiint \alpha(z,\theta) w_1(  \gamma_{z,\theta}(t) ,\dot \gamma_{z,\theta}(t)) f(\gamma_{z,\theta}(t))\bar\phi(z,\theta)\, \d t\, \d S_z\, \d\theta\\
  &= \iint \alpha^\sharp( x,v) w_1(x,v) f(x)\bar\phi^\sharp (x,v) J^\flat(x,v) \,\d x\, \d v,
\end{split}
\]
where $\alpha^\sharp(x,v) = \alpha(z(x,v), \theta(x,v))$, i.e., $\alpha^\sharp$ equals $\alpha$, extended as constant along the curves $\gamma_{z,\theta}(t)$; and the meaning of $\phi^\sharp$ is the same. Therefore,
\[
I_{\Gamma,w_1,\alpha}^* \phi(x) = \int_{|v|=1} \alpha^\sharp( x,v) \bar w_1(x,v) \phi^\sharp (x,v) J^\flat(x,v) \, \d  v.
\]
Let $w_0 \in C^\infty_0(T\Mint_1)$ be another function with small enough support. Then 
\be{43'}
\begin{split}
I_{\Gamma,w_1,\alpha}^* I_{\Gamma,w_0,\alpha}f(x)  = \int _{S^{n-1}}\int & (\alpha^\sharp\bar w_1)(x,v) 
(\alpha^\sharp w_0)\big(  \exp_x(t,v  )   ,  \partial_t\exp_x(t,v  ) \big)\\ &
\times f\big(\exp_x(t,v  )\big)
J^\flat(x,v) \,\d t\, \d v.
\end{split}
\ee

The simplicity assumption implies that for any $t_0\not=0$, and $(x_0,v_0)$ belonging to the support of the integrand above, the map $(t,v) \mapsto y=\exp_x(t,v)$ is a diffeomorphism from a neighborhood of $(t_0,v_0)$ to its image, and this is true for $x$ in some neighborhood of $x_0$. On the other hand, near $t_0=0$, the map $(t,v) \mapsto y=\exp_x(t,v)$ has Jacobian vanishing at $t=0$, and the ``true exponential map'' $\xi = tv  \mapsto y=\exp_x(t,v)$ is only a $C^1$ diffeomorphism, in general. 
By a compactness argument, given $\eps>0$, one can cover $M$ with finitely many charts, so that when $x$ belongs to either one of them, one can split the integration in \r{43'} into finitely many open sets that cover $\R\setminus(-\eps,\eps) \times S^{n-1}$. In each of those integrals, 
perform the change of variables $(t,v) \mapsto y=\exp_x(t,v)$. 
Then we get that the l.h.s.\ of \r{43'} is an operator with a smooth kernel. The only contribution to the singularities of the kernel may therefore come from $t\in (-\eps,\eps)$. 

To analyze the contribution to \r{43'} from $t\in(-\eps,\eps)$, we proceed as follows, see also \cite{DPSU}. The function 
$m(t,v;x)= (\exp_x(t,v) -x )/t$ is smooth, therefore
\begin{equation}   \label{S16}
\exp_x(t,v)-x =t m(t,v;x), \quad m(0,v;x) = \lambda(x,v) v.
\end{equation}
We introduce the new variables $(r,\omega) \in \R\times S^{n-1}$  by
\begin{equation}   \label{S17}
r = t|m(t,v;x)|, \quad \omega = m(t,v;x)/|m(t,v;x)|.
\end{equation}
Then $(r,\omega)$ are polar coordinates for $y-x =r\omega$ 
in which we allow $r$  to be ne\-gative. 
Clearly,  $(r,\omega)$ are smooth  at least for 
$\varepsilon$ small enough. 
Consider the Jacobian of this change of variables
\be{41'}
J(x,t,v) := \det \frac{\partial(r,\omega)}{\partial(t,v)},
\ee
computed with the same choice of local coordinates on $S^{n-1}$ for $v$ and $\omega$ (and independent of that choice). 
It is not hard to see that $J|_{t=0}= \lambda\not=0$, therefore the map 
$\R\times S^{n-1} \ni (t,v) \mapsto (r,\omega)\in \R\times S^{n-1} $ 
is a local  diffeomorphism from $\n(0)\times S^{n-1}$ to its image. 
We can decrease $\varepsilon$ if needed to ensure that it is 
a (global) diffeomorphism on its domain because then it is clearly injective.  
We denote the inverse functions by $t=t(x,r,\omega)$, $v=v(x,r,\omega)$. 
Note that in the $(r,\omega)$ variables
\begin{equation}   \label{S18}
t=r/\lambda +O(|r|), \quad v = \omega+O(|r|),\quad \dot \gamma_{x,v}(t) 
= \lambda\omega +O(|r|).
\end{equation}

Another representation of the new coordinates can be given by 
\[
r=\text{sign}(t)\left| \exp_x(t,v)  -x \right|, 
\quad 
\omega = \text{sign}(t)\frac{\exp_x(t,v)  -x}
{\left| \exp_x(t,v)  -x \right|},
\]
and
\[
(t,v) = \exp_x ^{-1}(x+r\omega)
\]
with the additional condition that $r$ and $t$ have the same sign 
(or are both zero). 

We return to \r{43'}. The paragraph after it shows that one can multiply the integrand by a smooth function $\chi(t)$ so that $\chi=1$ near $t=0$ and $\supp\chi$ is small enough; and the error is a smoothing operator. Then one can write, modulo a smoothing operator applied to $f$: 
\be{44}
\begin{split} 
I_{\Gamma,w_1,\alpha}^* I_{\Gamma,w_0,\alpha}f(x)  &\equiv \int _{S^{n-1}}\int \chi(t) B(x,t,v) f\big(\exp_x(t,v)\big) \, \d t \, \d v\\
&= \int_{S^{n-1}}\int \chi(t)  J^{-1}(x,t,v)B(x,t, v)f(x+r\omega)\big|_{t = t(x,r,\omega), \, v =v(x,r,\omega)} \, \d r\, \d\omega,
\end{split}
\ee
where, see \r{43'},
\be{45'}
B(x,t,v) = (\alpha^\sharp\bar w_1)(x,v)  (\alpha^\sharp w_0)\big(  \exp_x(t,v  )   ,  \partial_t\exp_x(t,v  ) \big) J^\flat(x,v).
\ee

\subsection{Certain class of integral operators 
with singular kernels}\label{singular} 

Let $U\subset \R^n$ be open and bounded. The integral representation \r{44} shows that we need to study integral operators with singular Schwartz kernels (with integrable singularity at $x=y$) of the class below, see also \cite[Appendix~D]{DPSU}.

\begin{lemma}   \label{lemma_SA}
Let $\mathcal{A} :C_0(U)\to C(U)$ be the operator 
\begin{equation}  \label{S19a}
\mathcal{A}f(x) = \int_{S^{n-1}}\int_{\R}  A(x,r,\omega) f(x+r\omega) \,\d r\,\d \omega,
\end{equation}
with $A\in C^\infty(U\times \R\times S^{n-1})$. 
Then $\mathcal A$ is a classical $\Psi$DO of order $-1$ with full symbol 
\[
a(x,\xi) \sim \sum_{k=0}^\infty a_k(x,\xi),
\]
where
\[
a_k(x,\xi) = 2\pi \frac{\i^k}{k!}\int_{S^{n-1}} \partial_r^k A(x,0,\omega)\delta^{(k)}(\omega\cdot \xi) \, \d\omega.
\]
\end{lemma}

\begin{proof} 
Notice first that if $A$ is an odd function of $(r,\omega)$, then $\mathcal Af=0$. 
Therefore, we can replace $A$ above by 
$A_{\text{even}}(r,\omega) = (A(r,\omega) + A(-r,-\omega))/2$. 
Next, it is easy to check that we can integrate over $r\ge0$ only and 
double the result. Therefore,
\begin{equation}  \label{S19b}
\mathcal Af(x) = 2\int_{S^{n-1}}\int_0^\infty  A_{\text{even}}(x,r,\omega) f(x+r\omega) \,\d r\,\d \omega.
\end{equation} 
Consider now $r$, $\omega$ as polar coordinates for $z=r\omega$, 
and make also the change of variables $y=x+z$ to get
\begin{equation}  \label{S20}
\mathcal Af(x) = 2 \int A_{\text{even}}\bigg(x,|y-x|  ,\frac{y-x}{|y-x|} \bigg) \frac{f(y)}{  |y-x|^{n-1}}\,\d y.
\end{equation}
Let
\begin{equation}  \label{S20a}
A_{\text{even}}(x,r,\omega) = \sum_{k=0}^{N-1} A_{\text{even},k} (x,\omega) r^k+r^{N}R_N(x,r,\omega)
\end{equation} 
be a finite Taylor expansion of $A_{\text{even}}$ in $r$ near $r=0$ with $N>0$. 
It follows easily that 
$2A_{\text{even},k}(x,\omega) = A_k(x,\omega) +(-1)^k A_k(x,-\omega)$, 
where $k!A_k = \partial^k_r|_{r=0}A$, and in particular, 
$A_{\text{even},k}(x,\omega)r^k$ is even w.r.t.\ $(r,\omega)$.  
The remainder term contributes to (\ref{S20}) an operator that 
maps $L^2_{\text{comp}}(U)$ into $H^{N-N_0}(U)$   with some fixed $N_0$. 
To study the contribution of the other terms, write
\begin{equation}  \label{S20b}
\mathcal A_{\text{even},k} f(x) = 2\int A_{\text{even},k}\bigg(x,  \frac{y-x}{|y-x|} \bigg) |y-x|^{k-n+1} f(y) \, \d y.
\end{equation}
The kernel of $\mathcal A_{\text{even},k}$ is therefore a function of $x$ and $z=y-x$, with a polynomial 
singularity at $y-x=0$, and it is therefore a formal $\Psi$DO with symbol 
that can be obtained by taking Fourier transform in the $z$ variable. 
Motivated by this, apply the Plancherel theorem to the integral above to get
\[
\mathcal A_{\text{even},k}f(x) = (2\pi)^{-n} \int e^{ix\cdot \xi} a_k(x,\xi) \hat f(\xi)\, \d\xi,
\]
where
\begin{align}  \nonumber
a_k(x,\xi) &= 2\int e^{-\i y\cdot \xi} A_{\text{even},k}\bigg(x, \frac{y-x}{|y-x|} \bigg) |y-x|^{k-n+1}   \,\d y\\ \nonumber
           &= 2\int_{S^{n-1}}\int_0^\infty e^{-\i r\omega\cdot \xi} A_{\text{even},k} (x,\omega) r^k\, \d r\, \d\omega\\ \nonumber
           &= \int_{S^{n-1}}\int_{-\infty}^\infty e^{-\i r\omega\cdot \xi} A_k (x,\omega) r^k\, \d r\, \d\omega\\
           &= 2\pi \i^k\int_{S^{n-1}}  A_k(x,\omega)\delta^{(k)}(\omega\cdot \xi) \, \d\omega.     \label{S21}
\end{align}
In the third line, we used the fact that 
$A_{\text{even},k}(x,\omega) r^k$ is even. 
Note that $a_k(x,\xi)$ is homogeneous in $\xi$ of order $-k-1$ 
and smooth away from $\xi=0$ but a distribution (in $\mathcal{S}'$) near zero. 
To deal with this, choose $\chi\in C_0^\infty$  supported in $|\xi|\le1$ 
and equal to $1$ near $\xi=0$. 
Write $a(x,\xi) = \chi(\xi)a(x,\xi) + (1-\chi(\xi))a(x,\xi)$. 
The second term is a classical amplitude, 
while the first one contributes the term
\begin{equation}  \label{S21a}
\mathcal A_{\text{even},k} (\check{\chi}*f )
\end{equation}
to \eqref{S20b} that is smooth,  as can be easily seen by making 
the change of variables $z=y-x$ in \eqref{S20b}. 
\end{proof}

\begin{proof}[Proof of Proposition~\ref{pr_2}] 
For $x$ in a small enough neighborhood of a fixed $x_0$, using a partition of unity $\{\chi_j\}$, we can express $N_{\Gamma,w,\alpha}$ as a finite sum of operators of the kind \r{43'}, namely $N_{\Gamma,w,\alpha} =\sum I^*_{\Gamma,w_j,\alpha} I^{}_{\Gamma,w_i,\alpha}$ with $w_i=\chi_i w$. By the analysis following \r{43'}, their Schwartz kernels are smooth if we integrate outside any interval containing $t=0$, and the only non-smooth contribution may come from terms of the kind \r{44}, where $w_0$ and $w_1$ are replaced by some $w_i$ and $w_j$.  By Lemma~\ref{lemma_SA}, \r{43'} is a classical  \PDO\ of order $-1$. Its principal symbol is given by 
\[
a_0(x,\xi) = 2\pi \int_{S^{n-1}} A(x,0,\omega)\delta(\omega\cdot \xi) \, \d\omega.
\]
In case of \r{44}, $A(x,0,\omega)$ is given by $A_{ij} =  J^{-1}(x,0,\omega)\big( \big|\alpha^\sharp w\big|^2 \chi_i \chi_j\big)(x,\omega)$. Then  $\sum_{ij} A_{ij}$ is elliptic because $w\not=0$, and because given $(x,\xi)$, there is $\omega\perp \xi$ so that $\alpha^\sharp(x,\omega)\not=0$, and there exists  $i$ so that $\chi_i(x,\omega)>0$; and all other terms are non-negative. Therefore, $N_{\Gamma,w,\alpha}$ is an elliptic \PDO\ of order $-1$ in $\Mint_1$, and the proposition follows. 
\end{proof}

The next proposition is a standard consequence of the ellipticity of $N_{\Gamma,w,\alpha}$. See \cite[Theorem~2]{SU-Duke} for a similar statement in tensor tomography. In contrast to \cite{SU-Duke} however, we do not lose a derivative.

\begin{proposition}   \label{pr_3}
Under the conditions of Theorem~\ref{thm_2}, without assuming that $I_{\Gamma,w,\alpha}$ is injective,  

(a)  one has the a~priori estimate
\[
\|f\|_{L^2(M)}  \le C
\|N_{\Gamma,w,\alpha}   f\|_{ H^1(M_1)} + C_s\|f\|_{H^{-s}(M_1)}, \quad \forall s;
\]

(b) $\text{\rm Ker}\, I_{\Gamma,w,\alpha}$ is finite dimensional and included in $C^\infty(M)$. 
\end{proposition}

\begin{proof}
Part (a) follows directly from Proposition~\ref{pr_2}. Next, if $f\in \text{\rm Ker}\, I_{\Gamma,w,\alpha}$, then $(\Id+K)f=0$, and $K$ is a compact operator on $L^2(M)$, with smooth kernel. This proves (b). 
\end{proof}

\section{Reducing the smoothness requirements}  
In this section, we will reduce the smoothness requirements on $\Gamma$ and the weight $w$, and will prove Theorem~\ref{thm_2}. 

We start with the observation that assuming that $I_{\Gamma,w,\alpha}$ is injective on $L^2(M)$, then $N_{\Gamma,w,\alpha}: L^2(M) \to H^1(M_1)$ is injective, also, see \cite{SU-Duke}. Then we  get by Proposition~\ref{pr_3}(b) and  \cite[Proposition~V.3.1]{Ta1} that  
\be{51}
\|f\|_{L^2(M)}  \le C
\|N_{\Gamma,w,\alpha}   f\|_{ H^1(M_1)}.
\ee
The second inequality in \r{est} is obvious. This proves part (a) of Theorem~\ref{thm_2}. 

In the rest of this section, we will perturb $\Gamma,w,\alpha$ and show that this will result in a small constant times $\|f\|_{L^2(M)}$ that can be absorbed by the l.h.s.\ above. We think of $\Gamma$ as determined by $(G,\mu,\sigma)$. Since $N_{\Gamma,w,\alpha}$ is a \PDO\ that depends on $\Gamma,w,\alpha$ in a continuous way, if the latter belongs to $C^k$, $k\gg2$, the statement of Theorem~\ref{thm_2}(b) follows immediately from what we already proved if $C^2$ there is replaced by $C^k$, $k\gg2$, see also \cite{SU-Duke, SU-rig, SU-AJM}. Our goal here is to reduce that smoothness requirement.

\begin{proposition}  \label{pr_4}
Assume that $G,\mu,\sigma,w,\alpha$ are fixed and belong to $C^2$. Let $(\tilde G,\tilde \mu,\tilde \sigma,\tilde w,\tilde \alpha)$ be $O(\delta)$ close to $(G,\mu,\sigma,w,\alpha)$ in $C^2$. Then there exists a constant $C>0$ that depends on an a priori bound on the $C^2$ norm of $(G,\mu,\sigma,w,\alpha)$, so that 
\be{51a}
\big\|\big(N_{\tilde \Gamma,\tilde w,\tilde \alpha} - N_{\Gamma,w,\alpha} \big)  f\big\|_{ H^1(M_1)} \le C\delta \|f\|_{L^2(M)}
\ee
\end{proposition}

\begin{proof}

Assume now that we have two systems  $(\tilde G,\tilde \mu,\tilde \sigma,\tilde w,\tilde \alpha)$ and  $(G,\mu,\sigma,w,\alpha)$, as in the proposition. Let $C_0$ be a bound on the $C^2$ norm of the first system. All constants below will depend on $C_0$. Let $\delta$ be as in the proposition.  To estimate the difference of those quantities related to the two systems, we will need the following comparison inequality for ODEs of Gronwall type.

\begin{lemma}  \label{ODE}
Let $x$, $\tilde x$ solve the ODE systems
\[
x' = F(t,x), \quad \tilde x' = \tilde F(t, \tilde x),
\]
where $F$, $\tilde F$ are continuous functions from $[0,T]\times U$ to a Banach space $\mathcal{B}$, where $U\subset \mathcal{B}$ is open. Let $F$ be Lipschitz w.r.t.\ $x$ with a Lipschitz constant $k>0$. Assume that
\[
\|F(t,x) - \tilde F(t,x)\|\le \delta, \quad \forall t\in [0,T], \;\forall x\in U.
\]
Assume that $x(t)$, $\tilde x(t)$ stay in $U$ for $0\le t\le T$. Then for $0\le t\le T$,
\be{52a}
\|x(t) -\tilde x(t)\|\le e^{kt}\|x(0)-\tilde x(0)\| + \frac{\delta}{k}\left( e^{kt}-1\right).
\ee
\end{lemma}

For a proof see \cite{CL}. Note that the lemma can be used to compare the derivatives of $x$ and $\tilde x$ w.r.t.\ the initial conditions by differentiating w.r.t.\ the initial conditions first, and then applying the lemma. Since the curves $\gamma$ solve the  equation \r{11} considered in the phase space, see \r{12},  we get under the assumptions of Proposition~\ref{pr_4},
\be{52a1}
\|\gamma_{x,v} - \tilde\gamma_{x,v}\|_{C^2} + \|\dot \gamma_{x,v} - \dot{\tilde{\gamma}}_{x,v}\|_{C^2}  \le C\delta,
\ee
where the $C^2$ norm is w.r.t.\ $(x,v,t)$. The inclusion of $t$ can be easily deduced from  equation \r{11}. 

Assume that $G,\mu,\sigma,w,\alpha$ are fixed and belong to $C^2$. We will determine first the smoothness of the functions $r(x,t,v)$ and $\omega(x,t,v)$ defined in \r{S17}. Since 
\be{52}
m(t,v;x) = \int_0^1 \dot\gamma_{x,v}(st)\, \d s, 
\ee
we get that $m$ and  $\dot m$ are $C^2$ functions of their arguments. By \r{S17}, we get that $\partial^j_t r$ and $\partial^j_t\omega$, $j=0,1$, are $C^2$, also. 
In particular, the inverse functions $t(x,r,\omega)$, $v(x,r,\omega)$ are $C^2$, too. On the other hand, $J$ and $J^\flat$ in \r{44}, \r{45'} are $C^1$. Moreover, the difference of those functions for the two systems is $O(\delta)$ in the corresponding norms. 

Let us analyze first \r{43'} in the case when the kernel there is multiplied by $1-\chi(t)$, compare with \r{44}. As explained in the paragraph following \r{43'}, we perform the change of variables $(t,v) \mapsto y=\exp_x(t,v)$ that is $C^2$ in our case, and its Jacobian is $C^1$. Moreover, the Jacobians for the two systems differ by $O(\delta)$ in the $C^1$ norm. 
Then we get an integral operator with a $C^1$ kernel, vanishing near the diagonal $x=y$. Clearly, such an operator maps $L^2(M)$ into $C^1(M_1)$. Moreover, the difference of two such operators, related to  $(\tilde G,\tilde \mu,\tilde \sigma,\tilde w,\tilde \alpha)$ and  $(G,\mu,\sigma,w,\alpha)$, respectively, has a norm $O(\delta)$. 

The more interesting case is what happens near the diagonal. To analyze this, we expand $B$ in \r{45'} and $J$, see \r{41'}  as 
\be{53}
B(x,t,v) = B_0(x,v) + tB_1(x,t,v),\quad J^{-1}(x,t,v) = J_0(x,v) + tJ_1(x,t,v).
\ee
The explicit expressions for $B_0, B_1, J_0, J_1$ are listed below:
\begin{align*} \label{54a}
B_0(x,v)& = \big((\alpha^\sharp)^2 \bar w_1 w_0\big)(x,v) J^\flat(x,v),\\
B_1(x,t,v)& = \int_0^1b_1(x,st,v) \d s, \quad \text{where} \\
b_1(x,t,v) &= \frac{\partial}{\partial t} (\alpha^\sharp\bar w_1)(x,v)  (\alpha^\sharp w_0)\big(  \exp_x(t,v  )   ,  (\partial_t\exp_x)(t,v  ) \big) J^\flat(x,v),\\
J_0(x,v) & = J^{-1}(x,0,v) = \lambda^{-1}(x,v),\\
J_1(x,t,v)& = \int_0^1 
\frac{\partial J^{-1}}{\partial t} (x,st,v)\, \d s.
\end{align*}  
Notice that  $B_0, B_1, J_0, J_1 \in C^1$. Moreover,  $\tilde B_0,\tilde  B_1,\tilde  J_0, \tilde J_1$ differ by them by $O(\delta)$ in the $C^1$ norm. 

Then for $A$ in \r{S19a}, we get
\be{54}
\begin{split}
A(x,r,\omega) & =\chi(t(x,r,\omega) J^{-1}(x,t(x,r,\omega),v(x,r,\omega) )  B(x,t(x,r,\omega) ,v(x,r,\omega))  \\
& =: A_0(x,\omega) +rA_1(x,r,\omega) .
\end{split}
\ee
Here $A_0(x,\omega)$ and $A_1(x,r,\omega)$ are $C^1$ functions of all variables, and we have used the fact that $t(x,r,\omega)/r$ is $C^1$, too. As above, we get 
\be{54b}
\|A_0(x,\omega)- \tilde A_0(x,\omega)\|_{C^1} + \|A_1(x,r, \omega)- \tilde A_1(x,r, \omega)\|_{C^1} \le C\delta.
\ee

Let $\mathcal{A}_0$, $\mathcal{A}_1$ be as in Lemma~\ref{lemma_SA} related to $A_1$ and $rA_1$, respectively. Then the Schwartz kernel of $\mathcal{A}_0$, see \r{S20}, is $  2{A}_{0,\text{even}}(x,\omega)r^{-n+1}$, where we use the notation
\[
r=|x-y|, \quad \omega=(y-x)/r.
\]
Therefore $  2{A}_{0,\text{even}}(x,\omega)r^{-n+1}$ has singularity of the type $r^{-n+1}$, while the kernel of  $\mathcal{A}_{1}$  has singularity of the type $r^{-n+2}$. To estimate the $H^1$ norm of $\mathcal{A}$, we need to analyze the operator with kernel $\partial_x \mathcal{A}$. We get that formally, $\partial_x \mathcal{A}_{0}$ is an operator with a non-integrable singularity of the type $r^{-n}$, while $\partial_x  {A}_{1,\text{even}}$ is an operator with kernel that still has an integrable singularity. Let now $\tilde{\mathcal{A}}_{0,1}$ be related to $(\tilde G,\tilde \mu,\tilde \sigma,\tilde w,\tilde \alpha)$.  The contribution of $\tilde{\mathcal{A}}_{1} -  \mathcal{A}_1  $ to \r{51a} is easy to estimate using \r{54b}. The remaining question is whether
\be{55}
\big\|\big( \tilde{\mathcal{A}}_0   -  \mathcal{A}_0 \big)  f\big\|_{ H^1(U)} \le C\delta \|f\|_{L^2(U)}, \quad f\in L^2(U), 
\ee
where $U$ is as in Lemma~\ref{lemma_SA}, and in our case, is a small enough open set in a fixed coordinate chart of $M_1$. 
We showed above that  $\mathcal{A}_0$ an operator with a weakly singular kernel, and $\partial_x\mathcal{A}_0$ is formally an operator with singular kernel. The continuity properties of the latter class are well studied, see e.g., \cite{Stein, MP}, and the integration is understood in principle value sense. By the Calder\'on-Zygmund Theorem, see, e.g., \cite[Theorem~X1.3.1]{MP}, \cite{Stein}, a singular operator with kernel $K(x,y) = \Omega(x,\omega)r^{-n}$ is bounded on $L^2(U)$, if $\Omega$ has a mean value 0 in the $\omega$ variable, and belongs to $L^\infty(U_x; \; L^2(S^{n-1}))$. Then the norm of that operator is bounded by $C\|\Omega\|$, where the latter norm is in $L^\infty(U_x; \; L^2(S^{n-1}))$.

In our case, we start with an operator with weakly singular kernel  $2{A}_{0,\text{even}}(x,\omega)r^{-n+1}$ that is even w.r.t.\ $\omega$, since ${A}_{0,\text{even}}$ it is independent of $r$. Therefore the $x$-derivative, if we differentiate the occurrence of $x$ in $r$ and $\omega$ only,  is an odd function of $\omega$. This makes the kernel $\partial_x (2{A}_{0,\text{even}}(x,\omega)r^{-n+1})$ a singular odd one, up to a weakly singular kernel. Now \cite[Theorem~XI.11.1]{MP} says that this is actually the kernel of $\partial_x \mathcal{A}_0$, and  by the Calder\'on-Zygmund Theorem, its $L^2\to L^2$ norm is bounded by $C\|A_{0}\|_{C^1}$. We apply now those arguments to 
$\tilde{\mathcal{A}}_0   -  \mathcal{A}_0$ with the aid of \r{54b}. This yields \r{55} and completes the proof of the proposition.
\end{proof}

\begin{proof}[Proof of Theorem~\ref{thm_2}] 
We already proved part (a) in \r{51}. 
Combine estimate \r{51} and Proposition~\ref{pr_4} to get
\[
\begin{split}
\|f\|_{L^2(M)}  & \le C
\|N_{\Gamma,w,\alpha}   f\|_{ H^1(M_1)}\\
&\le C\big\|N_{\tilde \Gamma,\tilde w,\tilde \alpha}  f\big\|_{ H^1(M_1)}+ 
C\big\|\big(N_{\tilde \Gamma,\tilde w,\tilde \alpha} - N_{\Gamma,w,\alpha} \big)  f\big\|_{ H^1(M_1)} \\
&\le C\big\|N_{\tilde \Gamma,\tilde w,\tilde \alpha}  f\big\|_{ H^1(M_1)}  + C\delta \|f\|_{L^2(M)}.
\end{split}
\]
This immediately implies Theorem~\ref{thm_2}(b). 
\end{proof}


\begin{thebibliography}{DPSU}
\frenchspacing

\bibitem[AR]{AR} {\sc Yu. Anikonov and V. Romanov}, On uniqueness of determination of a form of first degree by its integrals along geodesics, {\it J. Inv. Ill-Posed Problems}, {\bf 5}(1997), no. 6, 487--480. 


\bibitem[BG]{BG} {\sc I. N. Bernstein and M. L. Gerver}, Conditions on distinguishability of metrics by hodographs. {\em Methods and Algorithms of Interpretation of Seismological Information}, Computerized Seismology {\bf 13}, Nauka, Moscow, 50--73 (in Russian). 

\bibitem[B]{B} {\sc J. Boman}, An example of non-uniqueness for a generalized Radon transform, {\em  J. Anal. Math.}, {\bf 
61} (1993), 395--401.

\bibitem[BQ]{BQ} {\sc F. Boman and E. Quinto}, Support theorems for real-analytic Radon transforms, {\em Duke Math. J.} {\bf 55}(4)(1987), 943--948.

\bibitem[Ch]{Ch} {\sc E. Chappa}, On the characterization of the kernel of the geodesic X-ray transform, {\em Trans. Amer. Math. Soc.} {\bf 358}(2006), 4793-4807.

\bibitem[CL]{CL} {\sc  E. Coddington and N. Levinson},   
Theory of ordinary differential equations,   
Malabar, Fla.,  R.E. Krieger, 1984 


\bibitem[D]{D} \textsc{N. Dairbekov}, Integral geometry problem for nontrapping manifolds,  {\it Inverse Problems}, {\bf 22} (2006), no. 2, 431--445.

\bibitem[DPSU]{DPSU} {\sc N. Dairbekov, G. Paternain, P. Stefanov and G. Uhlmann},  Boundary rigidity problem in the presence of~a~magnetic field, preprint. 


\bibitem[Gr]{Gr} {\sc M. Gromov}, Filling Riemannian manifolds, {\it J. Diff. Geometry} {\bf 18}(1983), no. 1, 1--148. 

\bibitem[GrU]{GrU} {\sc A. Greenleaf and G. Uhlmann}, Nonlocal inversion formulas for the
X-ray transform, {\it Duke Math. J.} {\bf 58}(1989), no.~1, 205--240.

\bibitem[Gu]{Gu} {\sc V. Guillemin}, On some results of Gel'fand in
integral geometry, in {\it Pseudodifferential operators and
applications} (Notre Dame, Ind., 1984), pp. 149--155. Amer. Math. Soc.,
Providence, RI, 1985.

\bibitem[GuS1]{GuS1}
{\sc V. Guillemin and S. Sternberg}, {Geometric Asymptotics},
Mathematical Surveys and Monographs, vol.~14, American Mathematical Society,
Providence, Rhode Island, 1977.

\bibitem[GuS2]{GuS2}
{\sc V. Guillemin and S. Sternberg}, Some problems in integral geometry
and some related problems in micro-local analysis, {\it Am. J. Math.},
{\bf 101}(1979), 915--955.


\bibitem[KSU]{KSU} {\sc C. Kenig, J. Sj\"ostrand and G. Uhlmann}, The Calder\'on Problem with partial data, to appear in {\em  Ann. Math.}

\bibitem[MP]{MP}  {\sc S. Mikhlin and S. Pr\"ossdorf}, Singular Integral Operators, Springer, 1986.


\bibitem[MN]{MN} {\sc C. Morrey and L. Nirenberg}, On the analyticity of the solutions of linear elliptic systems of partial differential equations, {\em Comm. Pure Appl. Math.} {\bf 10}(1957), 271--290.

\bibitem[Mu1]{Mu} {\sc R. Mukhometov}, Inverse kinematic problem of seismic on the plane, {\em  Math. Problems of Geophysics, Akad. Nauk. SSSR}, Sibirsk. Otdel., Vychisl. Tsentr, Novosibirsk, {\bf 6}(2), 243--252 (1975). (in Russian).

\bibitem[Mu2]{Mu2} \bysame
The reconstruction problem of a two-dimensional Riemannian metric, and integral geometry (Russian), {\it Dokl. Akad. Nauk SSSR} {\bf 232}(1977), no. 1, 32--35.


\bibitem[Mu3]{Mu3} \bysame 
On a problem of reconstructing Riemannian metrics, {\it Siberian Math. J.} {\bf 22}(1982), no. 3, 420--433. 

\bibitem[MuR]{Mu-R}
{\sc R. G. Mukhometov and V. G. Romanov}, On the problem of finding an isotropic Riemannian metric in
an $n$-dimensional space (Russian), {\it Dokl. Akad. Nauk SSSR} 
{\bf 243}(1978), no. 1, 41--44.

\bibitem[Pe]{Pe} \textsc{L. Pestov}, Questions of well-posedness of the ray
tomography problems, Sib. Nauch. Izd., Novosibirsk (2003), (Russian).


\bibitem[Q]{Q} {\sc E. Quinto}, 
Radon transforms satisfying the Bolker assumption, in:  Proceedings of conference ``Seventy-five Years of Radon Transforms,'' International Press Co. Ltd., Hong Kong, pp. 263--270, 1994.

\bibitem[R]{R} {\sc V. Romanov}, Integral geometry on geodesics of an isotropic Riemannian metric, {\em Soviet Math. Dokl.} {\bf 19}(4), 847--851.


\bibitem[Sh1]{Sh}
{\sc V. Sharafutdinov}, Integral geometry of tensor fields, VSP, Utrech,
the Netherlands, 1994.


\bibitem[Sh2]{Sh-sib} \bysame, An integral geometry problem in a nonconvex domain, {\it Siberian Math. J.} {\bf 43}(6)(2002), 1159--1168.

\bibitem[Sh3]{Sh-2d} \bysame, Variations of Dirichlet-to-Neumann map
and deformation boundary rigidity of simple 2-manifolds, preprint.

\bibitem[Sj]{Sj-Ast} {\sc J. Sj\"ostrand}, Singularit\'es analytiques microlocales, Ast\'erique {\bf 95}(1982), 1--166.



\bibitem[SU1]{SU0} {\sc P. Stefanov and G. Uhlmann}, Stability estimates for the hyperbolic
Dirichlet to Neumann map in anisotropic media, {\em J. Funct. Anal.} {\bf 154}(2)
(1998), 330--358.

\bibitem[SU2]{SU1} 
\bysame, 
Rigidity for metrics with the same lengths of geodesics, {\em Math. Res. Lett.} {\bf 5}(1998), 83--96. 

\bibitem[SU3]{SU-Duke} 
\bysame, 
Stability estimates for the X-ray transform of tensor fields and boundary rigidity,  {\em Duke Math. J.} {\bf 123}(2004), 445--467.

\bibitem[SU4]{SU-rig} 
\bysame, 
Boundary rigidity and stability for generic simple metrics,  {\em J. Amer. Math. Soc.}  {\bf 18}(2005), 975--1003.

\bibitem[SU5]{SU-AJM} 
\bysame, 
Integral geometry of tensor fields on a class of non-simple Riemannian manifolds, preprint.

\bibitem[St]{Stein} {\sc E. Stein}, Singular Integrals and Differentiability Properties of Functions, Princeton Univ. Press, Princeton, 1970.

\bibitem[Ta1]{Ta1} {\sc M. Taylor}, Pseudodifferential Operators. {\em Princeton Mathematical Series} {\bf 34}. Princeton University Press, Princeton, N.J., 1981.



\bibitem[Tre]{T} {\sc F. Treves}, Introduction to Pseudodifferential and Fourier Integral Operators, Vol. 1. Pseudodifferential Operators. The University Series in Mathematics,  Plenum Press, New York--London, 1980. 


\end{thebibliography}
\end{document}